\documentclass[12pt]{article}

\usepackage{amsmath,amssymb,amsthm}

\begin{document}

\newcommand{\End}{{\rm{End}\ts}}
\newcommand{\non}{\nonumber}
\newcommand{\wt}{\widetilde}
\newcommand{\wh}{\widehat}
\newcommand{\ot}{\otimes}
\newcommand{\la}{\lambda}
\newcommand{\al}{\alpha}
\newcommand{\be}{\beta}
\newcommand{\ga}{\gamma}
\newcommand{\de}{\delta^{}}
\newcommand{\om}{\omega^{}}
\newcommand{\hra}{\hookrightarrow}
\newcommand{\ve}{\varepsilon}
\newcommand{\ts}{\,}
\newcommand{\qin}{q^{-1}}
\newcommand{\tss}{\hspace{1pt}}
\newcommand{\U}{ {\rm U}}
\newcommand{\Y}{ {\rm Y}}
\newcommand{\C}{\mathbb{C}\tss}
\newcommand{\ZZ}{\mathbb{Z}}
\newcommand{\A}{\mathcal{A}}
\newcommand{\Z}{{\rm Z}}
\newcommand{\gl}{\mathfrak{gl}}
\newcommand{\oa}{\mathfrak{o}}
\newcommand{\spa}{\mathfrak{sp}}
\newcommand{\g}{\mathfrak{g}}
\newcommand{\ka}{\mathfrak{k}}
\newcommand{\p}{\mathfrak{p}}
\newcommand{\sll}{\mathfrak{sl}}
\newcommand{\agot}{\mathfrak{a}}
\newcommand{\qdet}{ {\rm qdet}\ts}
\newcommand{\sdet}{ {\rm sdet}\ts}
\newcommand{\sgn}{ {\rm sgn}\ts}
\newcommand{\ev}{ {\rm ev}}
\newcommand{\Sym}{\mathfrak S}

\newcommand{\scl}{\scriptstyle}
\newcommand{\lan}{\langle\ts}
\newcommand{\ran}{\ts\rangle}
\newcommand{\Tau}{ {\mathcal T}}
\newcommand{\Yq}{{\Y_q(\gl_n)}}
\newcommand{\Uq}{{\U_q(\gl_n)}}
\newcommand{\eva}{ {\rm ev}_a}
\newcommand{\h}{\mathfrak h}
\newcommand{\ti}{T}
\newcommand{\La}{\Lambda}
\newcommand{\xiL}{\xi^{}_{\La}}
\newcommand{\etaL}{\eta^{}_{\La}}
\newcommand{\etaLo}{\eta^{}_{\La^0}}

\renewcommand{\theequation}{\arabic{section}.\arabic{equation}}

\newtheorem{thm}{Theorem}[section]
\newtheorem{lem}[thm]{Lemma}
\newtheorem{prop}[thm]{Proposition}
\newtheorem{cor}[thm]{Corollary}

\theoremstyle{definition}
\newtheorem{defin}[thm]{Definition}
\newtheorem{example}[thm]{Example}

\theoremstyle{remark}
\newtheorem{remark}[thm]{Remark}

\newcommand{\bth}{\begin{thm}}
\renewcommand{\eth}{\end{thm}}
\newcommand{\bpr}{\begin{prop}}
\newcommand{\epr}{\end{prop}}
\newcommand{\ble}{\begin{lem}}
\newcommand{\ele}{\end{lem}}
\newcommand{\bco}{\begin{cor}}
\newcommand{\eco}{\end{cor}}
\newcommand{\bde}{\begin{defin}}
\newcommand{\ede}{\end{defin}}
\newcommand{\bex}{\begin{example}}
\newcommand{\eex}{\end{example}}
\newcommand{\bre}{\begin{remark}}
\newcommand{\ere}{\end{remark}}

\newcommand{\bal}{\begin{aligned}}
\newcommand{\eal}{\end{aligned}}
\newcommand{\beq}{\begin{equation}}
\newcommand{\eeq}{\end{equation}}
\newcommand{\ben}{\begin{equation*}}
\newcommand{\een}{\end{equation*}}

\newcommand{\bpf}{\begin{proof}}
\newcommand{\epf}{\end{proof}}

\def\beql#1{\begin{equation}\label{#1}}

\title{\Large\bf  On irreducibility of tensor products of evaluation
modules for the quantum affine algebra}
\author{{A. I. Molev, \quad V. N. Tolstoy\quad and\quad R. B. Zhang}}

\date{} 


\maketitle

\vspace{15 mm}

\begin{abstract}
Every irreducible finite-dimensional representation of the
quantized enveloping algebra $\U_q(\gl_n)$ can be extended
to the quantum affine algebra $\U_q(\wh\gl_n)$
via the evaluation homomorphism.
We give in explicit form the necessary and sufficient conditions for
irreducibility of tensor products of such evaluation modules.

\end{abstract}

\vspace{39 mm}

\noindent School of Mathematics and Statistics\newline
University of Sydney, NSW 2006, Australia\newline
alexm@maths.usyd.edu.au

\vspace{7 mm}

\noindent
Institute of Nuclear Physics\newline
Moscow State University, Moscow 119992, Russia\newline
tolstoy@nucl-th.sinp.msu.ru

\vspace{7 mm}

\noindent School of Mathematics and Statistics\newline
University of Sydney, NSW 2006, Australia
\newline
rzhang@maths.usyd.edu.au

\newpage

\section{Introduction}
\label{sec:int}
\setcounter{equation}{0}

Let $\g$ be a simple Lie algebra over $\C$. Finite-dimensional irreducible
representations of the corresponding quantum affine algebra $\U_q(\wh\g)$
were classified by Chari and Pressley~\cite{cp:qa}; see also
\cite[Chapter~12]{cp:gq}. The representations are parameterized by
$n$-tuples of monic polynomials $(P_1(u),\dots,P_n(u))$, where $n$
is the rank of $\g$. Moreover,
every such representation is isomorphic to a subquotient
of a tensor product of the form
\beql{tenfund}
L_{a_1}(\om_{i_1})\ot\cdots\ot L_{a_k}(\om_{i_k}),
\end{equation}
where $L_a(\om_i)$ denotes
the so-called fundamental representation of $\U_q(\wh\g)$ which
corresponds to the $n$-tuple of polynomials with $P_j(u)=1$
for all $j\ne i$ and $P_i(u)=u-a$, where $a\in\C$ and $\om_i$
is a fundamental weight.
In general, the structure of
the tensor product module \eqref{tenfund}
appears to be rather
complicated. Only recently, irreducibility conditions
for this module were found.
These conditions were first conjectured by
Akasaka and Kashiwara in \cite{ak:fd} and proved there
in the cases where $\wh\g$ is of type $A^{(1)}$ or $C^{(1)}$.
In some other cases the conjecture was proved in different ways
by Frenkel--Mukhin~\cite{fm:cq}, Varagnolo--Vasserot~\cite{vv:sm}
before the general conjecture was settled by Kashiwara~\cite{k:lz}.
This result was generalized by Chari~\cite{c:bg} who gave, in particular,
irreducibility conditions for tensor products of the representations
corresponding to $n$-tuples of polynomials $(P_1(u),\dots,P_n(u))$
such that $P_j(u)=1$ for all $j\ne i$ and the
roots of $P_i(u)$ form a `$q$-string'.

In the case where the Lie algebra $\g$ is of type $A$, the corresponding quantum
affine algebra admits a class of evaluation modules.
Namely, as shown by Jimbo~\cite{j:qu},
there exists a family of algebra homomorphisms
$\eva:\U_q(\wh\g)\to\U_q(\g)$ which allows one to extend
any finite-dimensional irreducible representation $L(\la)$ of $\U_q(\g)$
to a module $L_a(\la)$ over $\U_q(\wh\g)$, where $\la$ is the highest weight
of $L(\la)$ and $a\in\C$.
In particular, if $\la=\om_i$
is a fundamental weight, the evaluation module coincides
with the fundamental representation $L_a(\omega_i)$.
Furthermore, any finite-dimensional irreducible
representation of $\U_q(\wh\g)$ is isomorphic to
a subquotient of a tensor product module
\beql{teneval}
L_{a_1}(\la^{(1)})\ot\cdots\ot L_{a_k}(\la^{(k)});
\end{equation}
see e.g. \cite[Section~12.2]{cp:gq}. Irreducible
tensor products of the form \eqref{teneval} thus provide
an explicit realization
of a wider class of representations in comparison
with the modules \eqref{tenfund}. In fact, in the case
of $\g=\sll_2$ every type $1$ irreducible finite-dimensional representation of
$\U_q(\wh\sll_2)$ is isomorphic to a module of the form \eqref{teneval};
see \cite[Section~12.2]{cp:gq}.

In this paper, we prove necessary and sufficient conditions
for irreducibility of the tensor product \eqref{teneval}.
It is more convenient for us to work with the quantum affine algebra
$\U_q(\wh\gl_n)$ instead of $\U_q(\wh\sll_n)$. The results can be
easily reformulated for the latter algebra as well.
Our starting point is the important {\it binary property\/}
established by Nazarov and Tarasov~\cite{nt:it} with the use
of an observation made by Kitanine, Maillet and Terras~\cite{kmt:ff, mt:oq}.
Namely, the representation \eqref{teneval}
is irreducible if and only if for all $i<j$ the modules
$L_{a_i}(\la^{(i)})\ot L_{a_j}(\la^{(j)})$ are irreducible.
Therefore, we only need to find an irreducibility criterion
for the case of $k=2$ factors in \eqref{teneval}.
In fact, the binary property is proved in \cite{nt:it}
in the Yangian context. An explicit formulation for the quantum affine algebra
case can be found in Leclerc, Nazarov
and Thibon~\cite{lnt:ir}.
A general binary cyclicity property was established by
Chari~\cite{c:bg} for tensor products of arbitrary irreducible
finite-dimensional representations of $\U_q(\wh\g)$
with $\g$ of any type.
An irreducibility criterion for
the induction products of evaluation modules over the affine Hecke
algebras of type $A$ was given in \cite{lnt:ir}.
It implies an irreducibility criterion
for the $\U_q(\wh\gl_n)$-module
$L_{a}(\la)\ot L_{b}(\mu)$
where the highest weights satisfy some extra conditions:
assuming that $\la$ and $\mu$
are partitions (one may do this without loss
of generality), one should require that the sum of the
lengths of $\lambda$ and $\mu$ does not exceed $n$.
In the case where $\la$ and $\mu$ are multiples
of fundamental weights, the irreducibility
conditions are given by Chari~\cite{c:bg}.

The same kind of irreducibility questions can be posed for the Yangians $\Y(\g)$;
see e.g. \cite[Section~12.1]{cp:gq}. An irreducibility criterion
of the tensor product of two arbitrary
evaluation modules $L_a(\la)\ot L_b(\mu)$
over the Yangian $\Y(\gl_n)$ is given in \cite{m:ic}.
The conditions on $\la$ and $\mu$
essentially coincide with those of \cite{lnt:ir}.
Some particular cases of this criterion were also established
in \cite{nt:tp}.
It has been known as a ``folklore theorem"
that the finite-dimensional
representation theories for the Yangian  and
the quantum affine algebras are essentially the same. A rigorous result
in that direction was recently proved by Varagnolo~\cite{v:qv}.
It allows one to establish an irreducibility
criterion for the quantum affine algebras
by using the Yangian criterion of \cite{m:ic}.
However, the proof in \cite{v:qv} uses rather involved
geometric arguments. The aim of this paper is to give
an independent direct proof of the irreducibility
criterion appropriately modifying
the arguments of \cite{m:ic}. In particular,
this requires a development of a $q$-analog of the quantum minor
techniques employed in \cite{m:ic}. This provides a quantum minor
realization of the lowering operators for the quantized algebra
$\U_q(\gl_n)$ and allows a new derivation of the $q$-analog
of the Gelfand--Tsetlin formulas; cf.
Jimbo~\cite{j:qr}, Ueno, Takebayashi and  Shibukawa~\cite{uts:gz},
Nazarov and Tarasov~\cite{nt:yg}, Tolstoy~\cite{t:ep}.

\medskip

The second author would like to thank the School of Mathematics and Statistics,
the University of Sydney, for the warm hospitality during his visit.
All authors gratefully acknowledge financial
support from the Australian Research Council.

\section{Quantized algebras and evaluation modules}
\setcounter{equation}{0}

Fix a complex parameter $q$ which is nonzero and not a root of unity.
Following Jimbo~\cite{j:qu},
we introduce the $q$-analog $\Uq$ of the universal enveloping
algebra $\U(\gl_n)$ as an associative algebra generated by
the elements $t_1,\dots,t_n,t_1^{-1},\dots,t_n^{-1}$, $e_1,\dots,e_{n-1}$ and
$f_1,\dots,f_{n-1}$ with the defining relations
\beql{defUq}
\bal
t_it_j=t_jt_i, &\qquad t_it_i^{-1}=t_i^{-1}t_i=1, \\
t_ie_jt_i^{-1}=e_j \tss q^{\delta_{ij}-\delta_{i,j+1}},&\qquad
t_if_jt_i^{-1}=f_j \tss q^{-\delta_{ij}+\delta_{i,j+1}},\\
[e_i,f_j]=\frac{k_i-k_i^{-1}}{q-\qin}\ts\delta_{ij},&\qquad \text{with}\ \
k_i=t_it_{i+1}^{-1},\\
[e_i,e_j]=[f_i,f_j]=0&,\qquad\text{if}\ \ |i-j|>1,\\
[e_i,[e_{i\pm 1},e_i]_q]_q&=[f_i,[f_{i\pm 1},f_i]_q]_q=0,
\eal
\end{equation}
where we have used the notation $[a,b]_{\zeta}=ab-\zeta\ts ba$.
The $q$-analogs of the root vectors are defined inductively by
\begin{alignat}{2}
{}&e_{i,i+1}=e_i,\qquad {}&&e_{i+1,i}=f_i,
\non\\
\label{rootv}
{}&e_{ij}=[e_{ik},e_{kj}]_q,\qquad &&\text{for}\ \ i<k<j,\\
{}&e_{ij}=[e_{ik},e_{kj}]_{\qin},\qquad &&\text{for}\ \ i>k>j.
\non
\end{alignat}

We shall also use an $R$-matrix
presentation of the algebra $\U_q(\gl_n)$; see e.g.
\cite{j:qu} and \cite{rtf:ql}.
Consider the $R$-matrix
\beql{rmatrixc}
R=q\ts\sum_i E_{ii}\ot E_{ii}+\sum_{i\ne j} E_{ii}\ot E_{jj}+
(q-\qin)\sum_{i<j}E_{ij}\ot E_{ji}
\end{equation}
which is an element of $\End\C^n\ot \End\C^n$, where
the $E_{ij}$ denote the standard matrix units and the indices run over
the set $\{1,\dots,n\}$.
The quantized enveloping algebra $\U_q(\gl_n)$ is generated
by the elements $t_{ij}$ and $\bar t_{ij}$ with $1\leq i,j\leq n$
subject to the relations
\beql{defrel}
\bal
t_{ij}&=\bar t_{ji}=0, \qquad 1 \leq i<j\leq n,\\
t_{ii}\ts \bar t_{ii}&=\bar t_{ii}\ts t_{ii}=1,\qquad 1\leq i\leq n,\\
R\ts T_1T_2&=T_2T_1R,\qquad R\ts \overline T_1\overline T_2=
\overline T_2\overline T_1R,\qquad
R\ts \overline T_1T_2=T_2\overline T_1R.
\eal
\end{equation}
Here $T$ and $\overline T$ are the matrices
\beql{matrt}
T=\sum_{i,j}t_{ij}\ot E_{ij},\qquad \overline T=\sum_{i,j}
\overline t_{ij}\ot E_{ij},
\end{equation}
which are regarded as elements of the algebra $\U_q(\gl_n)\ot \End\C^n$.
Both sides of each of the $R$-matrix relations in \eqref{defrel}
are elements of $\U_q(\gl_n)\ot \End\C^n\ot \End\C^n$ and the subscripts
of $T$ and $\overline T$ indicate the copies of $\End\C^n$ where
$T$ or $\overline T$ acts; e.g. $T_1=T\ot 1$.
An isomorphism between the two presentations is given by the formulas
\beql{isomgln}
t_i\mapsto t_{ii},\qquad t^{-1}_i\mapsto \bar t_{ii},\qquad
e_i\mapsto -\frac{\bar t_{i,i+1}\ts t_{ii}}{q-\qin},\qquad
f_i\mapsto \frac{\bar t_{ii}\ts t_{i+1,i}}{q-\qin}.
\end{equation}
We shall identify the corresponding elements of $\U_q(\gl_n)$
via this isomorphism.

We now introduce the quantum affine
algebra $\U_q(\wh\gl_n)$ following \cite{rtf:ql}; see also
\cite{d:sr}, \cite{fm:ha}.
By definition, $\U_q(\wh\gl_n)$ has countably many
generators $t_{ij}^{(r)}$ and $\bar t_{ij}^{\ts(r)}$ where
$1\leq i,j\leq n$ and $r$ runs over nonnegative integers.
They are combined into the matrices
\beql{taff}
T(u)=\sum_{i,j=1}^n t_{ij}(u)\ot E_{ij},\qquad
\overline T(u)=\sum_{i,j=1}^n \bar t_{ij}(u)\ot E_{ij},
\end{equation}
where $t_{ij}(u)$ and $\bar t_{ij}(u)$ are formal series
in $u^{-1}$ and $u$, respectively:
\beql{expa}
t_{ij}(u)=\sum_{r=0}^{\infty}t_{ij}^{(r)}\ts u^{-r},\qquad
\bar t_{ij}(u)=\sum_{r=0}^{\infty}\bar t_{ij}^{\ts(r)}\ts u^{r}.
\end{equation}
The defining relations are
\beql{defrelaff}
\bal
t_{ij}^{(0)}&=\bar t_{ji}^{\ts(0)}=0, \qquad 1 \leq i<j\leq n,\\
t_{ii}^{(0)}\ts \bar t_{ii}^{\ts(0)}&=\bar t_{ii}^{\ts(0)}
\ts t_{ii}^{(0)}=1,\qquad 1\leq i\leq n,\\
R(u,v)\ts T_1(u)T_2(v)&=T_2(v)T_1(u)R(u,v),\\
R(u,v)\ts \overline T_1(u)\overline T_2(v)&=
\overline T_2(v)\overline T_1(u)R(u,v),\\
R(u,v)\ts \overline T_1(u)T_2(v)&=T_2(v)\overline T_1(u)R(u,v).
\eal
\end{equation}
Here $R(u,v)$ is
the trigonometric $R$-matrix given by
\beql{trRm}
\bal
R(u,v)={}&(u-v)\sum_{i\ne j}E_{ii}\ot E_{jj}+(\qin u-q\tss v)
\sum_{i}E_{ii}\ot E_{ii} \\
{}+ {}&(\qin-q)\tss u\tss\sum_{i> j}E_{ij}\ot
E_{ji}+ (\qin-q)\tss v\tss\sum_{i< j}E_{ij}\ot E_{ji},
\eal
\end{equation}
where the  subscripts
of $T(u)$ and $\overline T(u)$
are interpreted in the same way as in \eqref{defrel}.

\bre
One could consider the more general centrally
extended quantum affine algebra
$\U_{q,c}(\wh\gl_n)$ instead of $\U_q(\wh\gl_n)$.
It is well-known, however, that the action of $c$ is trivial
in finite-dimensional irreducible representations; see e.g.
\cite[Chapter~12]{cp:gq}.
\ere

The defining relations
can easily be rewritten in terms of the
generators. In particular, the relations between the $t_{ij}(u)$
take the form
\begin{multline}\label{defrelt}
(q^{-\de_{ik}}u-q^{\de_{ik}}v)\ts t_{ij}(u)\ts t_{kl}(v)
+(\qin-q)\tss (u\tss\de_{i>k}+v\tss\de_{i<k})\ts t_{kj}(u)\ts t_{il}(v)\\
=(q^{-\de_{jl}}u-q^{\de_{jl}}v)\ts t_{kl}(v)\ts t_{ij}(u)
+(\qin-q)\tss (u\tss\de_{j<l}+v\tss\de_{j>l}) \ts t_{kj}(v)\ts t_{il}(u),
\end{multline}
where $\de_{i>k}$ or $\de_{i<k}$ is $1$ if the inequality
for the subscripts holds
and $0$ otherwise.

A family of the {\it evaluation homomorphisms\/}
$\eva:\U_q(\wh\gl_n)\to \U_q(\gl_n)$ is
defined by
\beql{eval}
T(u)\mapsto T-a\ts\overline T\ts u^{-1},\qquad
\overline T(u)\mapsto \overline T-a^{-1}\ts T\ts u,
\end{equation}
where $a$ is a nonzero complex number.
Note that the $R$-matrix $R(u,v)$ satisfies
$R(c\tss u,c\tss v)=c\tss R(u,v)$
for any nonzero $c\in \C$. Therefore, the mapping
\beql{scal}
T(u)\mapsto T(c\tss u), \qquad \overline T(u)\mapsto \overline T(c\tss u)
\end{equation}
defines an automorphism of the algebra $\U_q(\wh\gl_n)$.
Clearly, the homomorphism $\eva$ is the composition
of such an automorphism with $c=a^{-1}$ and the evaluation
homomorphism $\ev=\ev_1$ given by
\beql{eval1}
T(u)\mapsto T-\overline T\ts u^{-1},\qquad
\overline T(u)\mapsto \overline T-\ts T\ts u.
\end{equation}

There is a Hopf algebra structure on $\U_q(\wh\gl_n)$ with the coproduct
defined by
\beql{copraff}
\Delta\big(t_{ij}(u)\big)=\sum_{k=1}^n t_{ik}(u)\ot t_{kj}(u),
\qquad
\Delta\big(\bar t_{ij}(u)\big)=\sum_{k=1}^n\bar t_{ik}(u)\ot\bar t_{kj}(u).
\end{equation}

Finite-dimensional irreducible representations of $\Uq$ are completely
described by their highest weights; see e.g. \cite[Chapter~10]{cp:gq}.
Let $\la=(\la_1,\dots,\la_n)$ be an $n$-tuple of integers with
the condition $\la_1\geq\cdots\geq\la_n$. We shall call such an $n$-tuple
a $\gl_n$-{\it highest weight\/}.
The finite-dimensional irreducible $\Uq$-module $L(\la)$
corresponding to the highest weight $\la$
contains a unique, up to a constant factor,
vector $\xi$ such that
\beql{hive}
t_i\ts \xi=q^{\la_i}
\ts \xi, \qquad\text{and}\qquad e_{ij}\ts \xi=0\quad\text{for}\ \ i<j.
\end{equation}
Using the evaluation homomorphism \eqref{eval} we can extend $L(\la)$
to the quantum affine algebra $\U_q(\wh\gl_n)$. We call such an extension
the {\it evaluation module\/} and denote it by $L_a(\la)$.
In the case $a=1$ we keep the notation $L(\la)$ for the
evaluation module $L_1(\la)$.
The coproduct \eqref{copraff} allows us to form
tensor product modules of the type $L_a(\la)\ot L_b(\mu)$
over the algebra $\U_q(\wh\gl_n)$.
Our main result is an irreducibility criterion of such modules.
We note that without loss of generality both
evaluation parameters $a$ and $b$ can be taken to be equal to $1$; see
Section~\ref{sec:emy} below. In order to formulate the result, we need
the following definition \cite{m:ic}; cf. \cite{lnt:ir}.
Two disjoint finite subsets $A$ and $B$ of $\ZZ$ are called {\it crossing\/}
if there exist elements $a_1,a_2\in A$ and $b_1,b_2\in B$ such that
either $a_1<b_1<a_2<b_2$, or $b_1<a_1<b_2<a_2$. Otherwise, $A$ and $B$
are called {\it non-crossing\/}.
Given a highest weight $\la$ with integer entries
we set $l_i=\la_i-i+1$ and
introduce the
following subset of $\ZZ$:
\beql{alambda}
\A_{\la}=\{l_1,l_2,\dots,l_n\}.
\non
\end{equation}

\bth\label{thm:main}
The $\U_q(\wh\gl_n)$-module $L(\la)\ot L(\mu)$ is irreducible if and only if
the sets $\A_{\la}\setminus \A_{\mu}$ and $\A_{\mu}\setminus \A_{\la}$
are non-crossing.
\eth

The proof of this theorem will be given in the next sections.
The first step is to reduce the problem to the case
of the $q$-{\it Yangian\/} $\Yq$. Then we develop
an appropriate $q$-version of the techniques
of lowering operators and Gelfand--Tsetlin bases
used in \cite{m:ic} for the proof of such criterion
in the Yangian case.

\section{Evaluation modules over $q$-Yangians}
\label{sec:emy}
\setcounter{equation}{0}

We define the
$q$-Yangian $\Yq$ as the (Hopf) subalgebra of $\U_q(\wh\gl_n)$
generated by the
elements $t_{ij}^{(r)}$ with $1\leq i,j\leq n$ and $r\geq 0$.
The restriction of the evaluation homomorphism \eqref{eval}
to the $q$-Yangian is given by the first formula in \eqref{eval},
or, equivalently, in terms of the first presentation of $\Uq$,
it can be written as
\begin{alignat}{2}
{}&t_{ii}(u)\mapsto t_i-a\tss t_i^{-1}u^{-1},&&
\non\\
\label{evalyang}
{}&t_{ij}(u)\mapsto (q-\qin)\ts t_j\tss e_{ij},\qquad &&\text{if}\ \ i>j,   \\
{}&t_{ij}(u)\mapsto a\tss (q-\qin)\ts e_{ij}\tss
t_i^{-1}u^{-1},\qquad &&\text{if}\ \ i<j.
\non
\end{alignat}

The highest weight of an arbitrary
finite-dimensional irreducible representation $L$ of $\Uq$
may have a more general form than \eqref{hive}.
Namely, if $\xi$ is the highest vector of $L$
then
\beql{hive2}
t_i\ts \xi=\al_i\ts \xi, \qquad\text{and}\qquad e_{ij}\ts
\xi=0\quad\text{for}\ \ i<j,
\end{equation}
for a collection $(\al_1,\dots,\al_n)$ of nonzero complex numbers
of the form
\beql{hvform}
\al_i=h\tss\ve_i\tss q^{\la_i},\qquad i=1,\dots,n,
\end{equation}
where $h$ is a nonzero complex number,
each $\ve_i$ is equal to $1$ or $-1$, and
the $\la_i$ are integers satisfying $\la_i\geq \la_{i+1}$ for all $i$.
We denote the corresponding representation by $L(h,\ve,\la)$ where
\beql{vela}
\ve=(\ve_1,\dots,\ve_n),\qquad \la=(\la_1,\dots,\la_n).
\end{equation}
The evaluation homomorphism $\eva$ allows us
to regard $L(h,\ve,\la)$
as a $\Yq$-module which we denote by $L_a(h,\ve,\la)$.
Using the coproduct \eqref{copraff}
we can consider the tensor products of the form
\beql{ten2}
L_a(h,\ve,\la)\ot L_{a'}(h',\ve',\la')
\end{equation}
as $\Yq$-modules. It is clear that this $\Yq$-module coincides
with the restriction of \eqref{ten2} regarded as a $\U_q(\wh\gl_n)$-module.

\bpr\label{prop:qatoyang}
The $\U_q(\wh\gl_n)$-module \eqref{ten2} is irreducible if and only if
its restriction to the $q$-Yangian $\Yq$ is irreducible.
\epr

\bpf
The ``if" part is obviously true. Now
suppose, on the contrary, that the $\Yq$-module \eqref{ten2}
contains a nontrivial submodule $W$. Then,
by \eqref{eval} and \eqref{copraff}, $W$ is invariant with respect to
all operators of the form
\beql{tenoper}
\sum_{k=1}^n (t_{ik}-\bar t_{ik}\ts u^{-1}\ts a)\ot
(t_{kj}-\bar t_{kj}\ts u^{-1}\ts {a'}),\qquad i,j=1,\dots,n.
\end{equation}
This implies that $W$ is invariant with respect to
the operators
\beql{tenoper2}
\sum_{k=1}^n (\bar t_{ik}-t_{ik}\ts u\ts a^{-1})\ot
(\bar t_{kj}-t_{kj}\ts u\ts {a'}^{-1}).
\end{equation}
However, the latter operator is the image of the
generator $\bar t_{ij}(u)$ of $\U_q(\wh\gl_n)$ in the module
\eqref{ten2}. Therefore, $W$ is invariant under the action of
the entire algebra
$\U_q(\wh\gl_n)$. But this contradicts the irreducibility
of \eqref{ten2}.
\epf

\bre
Since every finite-dimensional irreducible representation $V$ of
$\U_q(\wh\gl_n)$ is isomorphic to a subquotient of \eqref{teneval},
the above argument can obviously be extended to show that
any such $V$ remains irreducible when restricted to $\Y_q(\gl_n)$.
\ere

Proposition~\ref{prop:qatoyang} tells us that the irreducibility conditions
for tensor products of evaluation modules over $\U_q(\wh\gl_n)$
and $\Y_q(\gl_n)$ are the same. In what follows we work with
$\Y_q(\gl_n)$-modules.
Finite-dimensional irreducible
representations of $\Yq$ can be described in terms of their highest weights
in a way similar to the case of the Yangian $\Y(\gl_n)$
\cite{d:nr}; see also \cite[Chapter~12]{cp:gq} and \cite{m:fd}.
The highest weight of such a module $L$ is a collection
of formal power series $(\nu_1(u),\dots,\nu_n(u))$ in $u^{-1}$ such that
\beql{yanhv}
t_{ii}(u)\ts\zeta=\nu_i(u)\ts\zeta,\qquad\text{and}\qquad t_{ij}(u)\ts\zeta=0
\quad \text{for}\ \ i<j
\end{equation}
for a vector $\zeta\in L$ (the highest vector)
which is determined uniquely up to a constant factor.
If the $\Y(\gl_n)$-module \eqref{ten2} is irreducible
its highest weight is easy to find from
\eqref{copraff} and \eqref{evalyang}. It is given by
\beql{hvtpr}
\nu_i(u)=(\al_i-a\tss \al_i^{-1}\tss u^{-1})\ts
(\al'_i-a'\tss \al_i^{\prime -1}\tss u^{-1}),
\end{equation}
where the $\al_i$ and $\al_i'$ are the components of the highest weights
of $L_a(h,\ve,\la)$ and $L_{a'}(h',\ve',\la')$; see \eqref{hvform}.
Note that for a given nondegenerate diagonal
matrix $D={\rm diag\tss}(d_1,\dots,d_n)$,
the mapping
\beql{diag}
T(u)\mapsto D\ts T(u)
\end{equation}
defines an algebra automorphism of $\Yq$, as follows from
\eqref{defrelt}.
Taking the composition of \eqref{ten2} with this automorphism
where $d_i=h^{-1}h^{\prime -1}\ve_i\ve_i'$ we find that
the irreducibility of \eqref{ten2} is equivalent to that of the module
\beql{simte}
L_b(\la)\ot L_{b'}(\la'),\qquad b=a\tss h^{-2},\quad b'=a'\tss h^{\prime -2},
\end{equation}
where by $L(\la)$ we denote the $\Uq$-module $L(h,\ve,\la)$ with $h=1$ and
$\ve=(1,\dots,1)$. Similarly, using the automorphism \eqref{scal}
we find that for any nonzero $c\in \C$ the module \eqref{simte} is irreducible
if and only if
the module $L_{bc}(\la)\ot L_{b'c}(\la')$ is.  On the other hand, the module
\eqref{simte} is irreducible unless $b\tss b^{\prime -1}\in q^{2\ZZ}$.
Analogs of this fact are well-known
both in the case of Yangians and the quantum affine algebras; cf.
\cite[Chapter~12]{cp:gq}.
One of the ways to prove this is to show that both the module and its dual
have no nontrivial singular vectors by considering
the eigenvalues of the quantum determinant on the module;
see \eqref{qdet} below.
So, we may now assume that $b'=1$ and
$b=q^{2k}$ in \eqref{simte} for some $k\in \ZZ$.
However, using the automorphism \eqref{scal} with $d_i\equiv q^{-k}$
we conclude that the irreducibility of \eqref{simte} is equivalent to that
of the module $L_1(\la)\ot L_1(\la'-k I)$ where $I=(1,\dots,1)$.
We shall keep the notation $L(\la)$
for the $\Yq$-module $L_b(\la)$ with $b=1$.
Thus, we may only consider, without loss
of generality, the $\Yq$-modules
of the form $L(\la)\ot L(\mu)$. The irreducibility conditions
for a general tensor product module \eqref{ten2}
can be easily obtained from
Theorem~\ref{thm:main}.

\section{Quantum minor relations}
\label{sec:qmr}

Here we formulate some well-known properties of quantum determinants
and quantum minors; see e.g. \cite{c:ni}, \cite{nt:yg}.

Let us consider
the multiple tensor product
$\Yq\ot (\End\C^n)^{\ot\tss r}$.
We have the following corollary of \eqref{defrelaff}
which is verified in the
same way as for the Yangians; cf. \cite{mno:yc}:
\beql{fundam}
R(u_1,\dots,u_r)\ts T_1(u_1)\cdots T_r(u_r)=
T_r(u_r)\cdots T_1(u_1)\ts  R(u_1,\dots,u_r),
\end{equation}
where
\beql{Rlong}
R(u_1,\dots,u_r)=\prod_{i<j}R_{ij}(u_i,u_j),
\end{equation}
with the product taken in the lexicographical order on the pairs $(i,j)$.
Here, like in \eqref{defrelaff},
the subscripts of the matrices $T(u)$ and $R(u,v)$
indicate the copies of $\End\C^n$.
Consider the $q$-permutation operator $P\in\End(\C^n\ot\C^n)$
defined by
\beql{qperm}
P=\sum_{i}E_{ii}\ot E_{ii}+ q\tss\sum_{i> j}E_{ij}\ot
E_{ji}+ \qin\sum_{i< j}E_{ij}\ot E_{ji}.
\end{equation}
An action of the symmetric group $\Sym_r$ on the space $(\C^n)^{\ot\tss r}$
can be defined by setting $s_i\mapsto P_{s_i}:=P_{i,i+1}$ for $i=1,\dots,r-1$,
where $s_i$ denotes the transposition $(i,i+1)$.
If $\sigma=s_{i_1}\cdots s_{i_l}$ is a reduced decomposition
of an element $\sigma\in \Sym_r$ we set
$P_{\sigma}=P_{s_{i_1}}\cdots P_{s_{i_l}}$.
We denote by $A_r$ the $q$-antisymmetrizer
\beql{antisym}
A_r=\sum_{\sigma\in\Sym_r}\sgn\sigma\cdot P_{\sigma}.
\end{equation}
The following proposition is proved by induction on $r$ in the same way as for the
Yangians \cite{mno:yc} with the use of a property of the reduced decompositions
\cite{h:il}, p.50.

\bpr\label{prop:antisym}
We have the relation in $\ts\End (\C^n)^{\ot\tss r}$:
\beql{ranti}
R(1,q^{-2},\dots,q^{-2r+2})=\prod_{0\leq i<j\leq r-1}(q^{-2i}-q^{-2j})\ts A_r.
\end{equation}
\epr

Now \eqref{fundam} implies that
\beql{anitt}
A_r\ts T_1(u)\cdots T_r(q^{-2r+2}u)=T_r(q^{-2r+2}u)\cdots T_1(u)\ts  A_r
\end{equation}
which equals
\beql{matelmi}
\sum_{a_i,b_i}{t\ts}^{a_1\cdots\ts a_r}_{b_1\cdots\ts b_r}(u)\ot
E_{a_1b_1}\ot\cdots\ot E_{a_rb_r}
\end{equation}
for some elements ${t\ts}^{a_1\cdots\ts a_r}_{b_1\cdots\ts b_r}(u)\in\Yq[[u^{-1}]]$
which we call the {\it quantum minors\/}.
They can be given by the following formulas which are immediate from the definition.
If $a_1<\cdots<a_r$ then
\beql{qminor}
{t\ts}^{a_1\cdots\ts a_r}_{b_1\cdots\ts b_r}(u)=
\sum_{\sigma\in \Sym_r} (-q)^{-l(\sigma)} \cdot t_{a_{\sigma(1)}b_1}(u)\cdots
t_{a_{\sigma(r)}b_r}(q^{-2r+2}u),
\end{equation}
and for any $\tau\in\Sym_r$ we have
\beql{qmsym}
{t\ts}^{a_{\tau(1)}\cdots\ts a_{\tau(r)}}_{b_1\cdots\ts b_r}(u)=
(-q)^{l(\tau)}{t\ts}^{a_1\cdots\ts a_r}_{b_1\cdots\ts b_r}(u).
\end{equation}
Here $l(\sigma)$ denotes the length of the permutation $\sigma$.
If $b_1<\cdots<b_r$ (and the $a_i$ are arbitrary) then
\beql{qminor2}
{t\ts}^{a_1\cdots\ts a_r}_{b_1\cdots\ts b_r}(u)=
\sum_{\sigma\in \Sym_r} (-q)^{l(\sigma)} \cdot t_{a_rb_{\sigma(r)}}(q^{-2r+2}u)\cdots
t_{a_1b_{\sigma(1)}}(u),
\end{equation}
and for any $\tau\in\Sym_r$ we have
\beql{qmsym2}
{t\ts}^{a_1\cdots\ts a_r}_{b_{\tau(1)}\cdots\ts b_{\tau(r)}}(u)=
(-q)^{-l(\tau)}{t\ts}^{a_1\cdots\ts a_r}_{b_1\cdots\ts b_r}(u).
\end{equation}
Note also that the quantum minor is zero if two top or two bottom indices
are equal.

As another application of \eqref{fundam} we obtain
the relations between the generators $t_{ij}(u)$
and the quantum minors. For this we introduce an extra copy of
$\End\C^n$ as a tensor factor which will be enumerated by the index $0$.
Now specialize
the parameters $u_i$ as follows:
\beql{specpar}
u_0=u, \qquad u_i=q^{-2i+2}v\quad\text{for}\ \ i=1,\dots,r.
\end{equation}
Then by Proposition~\ref{prop:antisym}
the element \eqref{Rlong} will take the form
\beql{RA}
R(u,v,\dots,q^{-2r+2}v)=\prod_{i=1}^r R_{0i}(u,q^{-2i+2}v)\ts A_r.
\end{equation}
Using the definition of the quantum minors \eqref{matelmi}
and equating the matrix elements on both sides of \eqref{fundam}
we get the required relations.
To write them down, let us fix indices
$a,b$, $c_1<\cdots<c_r$ and $d_1<\cdots<d_r$. Then we have
\beql{aq=bq}
A_{a,b,(c),(d)}(u,v)=B_{a,b,(c),(d)}(u,v),
\end{equation}
where
\beql{alhs}
\bal
A_{a,b,(c),(d)}(u,v)&=(u-v)\ts t_{ab}(u)\ts
{t\ts}^{c_1\cdots\ts c_r}_{d_1\cdots\ts d_r}(v)\\
{}&+(\qin-q)\ts u\sum_{i=1}^k(-q)^{k-i}\ts t_{c_ib}(u)\ts
{t\ts}^{c_1\cdots\ts \wh{c}_i\cdots\ts  c_k a\ts
c_{k+1}\cdots\ts c_r}_{d_1\cdots\ts d_r}(v)\\
{}&+(\qin-q)\ts v\sum_{i=k+1}^r(-q)^{k-i+1}\ts t_{c_ib}(u) \ts
{t\ts}^{c_1\cdots\ts  c_k a\ts c_{k+1}\cdots\ts
\wh{c}_i\cdots \ts c_r}_{d_1\cdots\ts d_r}(v),
\eal
\end{equation}
if $c_k<a<c_{k+1}$ for some $k\in\{0,1,\dots,r\}$, and
\beql{alhs2}
A_{a,b,(c),(d)}(u,v)=(\qin u-q\tss v)\ts t_{ab}(u)\ts
{t\ts}^{c_1\cdots\ts c_r}_{d_1\cdots\ts d_r}(v),
\end{equation}
if $a=c_k$ for some $k$. Furthermore,
\beql{blhs}
\bal
B_{a,b,(c),(d)}(u,v)&=(u-v)\ts {t\ts}^{c_1\cdots\ts c_r}_{d_1\cdots\ts d_r}(v)
\ts t_{ab}(u)\\
{}&+(\qin-q)\ts v\sum_{i=1}^l(-q)^{i-l}\ts
{t\ts}^{c_1\cdots\ts c_r}_{d_1\cdots\ts \wh{d_i}\cdots
\ts d_l b\ts d_{l+1}\cdots\ts d_r}(v)
\ts t_{ad_l}(u)\\
{}&+(\qin-q)\ts u\sum_{i=l+1}^r(-q)^{i-l-1}\ts
{t\ts}^{c_1\cdots\ts c_r}_{d_1\cdots \ts d_l b\ts d_{l+1}
\cdots\ts \wh{d_i}\cdots \ts d_r}(v)
\ts t_{ad_l}(u),
\eal
\end{equation}
if $d_l<b<d_{l+1}$ for some $l\in\{0,1,\dots,r\}$, and
\beql{blhs2}
B_{a,b,(c),(d)}(u,v)=(\qin u-q\tss v)\ts
{t\ts}^{c_1\cdots\ts c_r}_{d_1\cdots\ts d_r}(v)
\ts t_{ab}(u),
\end{equation}
if $b=d_l$ for some $l$; the hats indicate that the indices
are omitted.
In particular, \eqref{aq=bq} implies the well-known
property of the quantum minors: for any indices $i,j$ we have
\beql{center}
[{t}^{}_{c_id_j}(u),
{t\ts}^{c_1\cdots\ts c_r}_{d_1\cdots\ts d_r}(v)]=0.
\end{equation}
This implies that all the coefficients of the series
\beql{qdet}
\qdet T(u)={t\ts}^{1\cdots\ts n}_{1\cdots\ts n}(u)
\end{equation}
belong to the center of $\Yq$. The element \eqref{qdet}
is called the {\it quantum determinant\/} of the matrix $T(u)$.
The {\it quantum comatrix\/} $\wh T(u)$
is defined by the relation
\beql{decom}
\wh T(u)\ts T(q^{-2n+2}u)=\qdet T(u).
\end{equation}
Using the definition \eqref{matelmi} we find that
$
\wh T(u)=\sum_{i,j}\wh t_{ij}(u)\ot E_{ij}
$
where
\beql{wtil}
\wh t_{ij}(u)=
(-q)^{j-i}\ts{t\ts}^{1\cdots\ts\wh j\ts\cdots n}_{1\cdots\ts\wh i\ts\cdots\ts n}(u).
\end{equation}
The elements $\wh t_{ij}(u)$ satisfy quadratic relations
which can be written in the $R$-matrix form
\beql{R21}
R(u,v)\ts \wh T_2(v)\ts  \wh T_1(u)=\wh T_1(u)\ts \wh T_2(v)\ts R(u,v).
\end{equation}
To see this, we multiply both sides of the first
$R$-matrix relation in \eqref{defrelaff}
by the product $T_1(u)^{-1}\ts T_2(v)^{-1}$
from  the left and by $T_2(v)^{-1}\ts T_1(u)^{-1}$ from the right.
Then substitute $u\mapsto q^{-2n+2}u$, $v\mapsto q^{-2n+2}v$
and multiply both sides by $\qdet T(u)\ts\qdet T(v)$ to get \eqref{R21}.
It is easy to rewrite this relation in terms of $\wh t_{ij}(u)$
in a way similar to \eqref{defrelt}.
We shall also need the following relation between the quantum minors.

\bpr\label{prop:1rel}
We have
\beql{1rel}
{t\ts}^{2\ts\cdots\ts n}_{1\ts\cdots\ts n-2,n}(u)\ts
{t\ts}^{2\ts\cdots\ts n-1}_{2\ts\cdots\ts n-1}(u)=
{t\ts}^{2\ts\cdots\ts n-1}_{1\ts\cdots\ts n-2}(u)
\ts{t\ts}^{2\ts\cdots\ts n}_{2\ts\cdots\ts n}(u)
+ q\ts{t\ts}^{2\ts\cdots\ts n}_{1\ts\cdots\ts n-1}(u)
\ts{t\ts}^{2\ts\cdots\ts n-1}_{2\ts\cdots\ts n-2,n}(u).
\end{equation}
\epr

\bpf We find from the definition of the quantum determinant that
\beql{dobtt}
A_n T_1(u)\cdots T_{n-2}(q^{-2n+6}u)=\qdet T(u)\ts
A_n T_{n}(q^{-2n+2}u)^{-1}T_{n-1}(q^{-2n+4}u)^{-1}.
\end{equation}
Equating the matrix elements of both sides and using \eqref{decom}
we arrive at the following relation: for $i<j$ and $k<l$,
\beql{tdobe}
(-q)^{i+j-k-l}\ts {t\ts}^{1\ts\cdots\ts \wh i\ts \cdots
\ts \wh j\ts \cdots\ts n}_{1\ts\cdots\ts \wh k\ts \cdots
\ts \wh l\ts \cdots\ts n}(u)\ts\qdet T(q^2u)
=\wh t_{lj}(u)\ts \wh t_{ki}(q^2u)-\qin \ts\wh t_{kj}(u)\ts \wh t_{li}(q^2u).
\end{equation}
The right hand side is
a $2\times 2$-quantum minor of the matrix $\wh T(u)$
and we denote it by ${\wh t\ }^{kl}_{ij}(u)$.
Furthermore, by analogy with \eqref{fundam} we obtain from \eqref{R21}
\beql{funh}
R(u_1,u_2,u_3)\ts \wh T_3(u_3)\ts \wh T_2(u_2)\ts  \wh T_1(u_1)=
\wh T_1(u_1)\ts \wh T_2(u_2)\ts \wh T_3(u_3)\ts  R(u_1,u_2,u_3).
\end{equation}
Now specialize $u_1=q^2v$, $u_2=v$, $u_3=u$ and take
the coefficients at $E_{n-1,1}\ot E_{nn}\ot E_{11}$ on both sides.
Using Proposition~\ref{prop:antisym}
and \eqref{tdobe}
we get
\beql{tiltt}
\bal
{}&(\qin-q)\ts v\ts {\wh t}_{n-1,1}(u)\ts
{\wh t\ }^{1n}_{1n}(v)-(1-q^2)\ts v\ts{\wh t}_{n1}(u)\ts
{\wh t\ }^{1,n-1}_{1n}(v)+(v-u)\ts{\wh t}_{11}(u)\ts
{\wh t\ }^{n-1,n}_{1n}(v)\\
{}={}&(\qin v-q\tss u)\ts {\wh t\ }^{n-1,n}_{1n}(v)\ts  {\wh t}_{11}(u).
\eal
\end{equation}
By putting $u=v$ and rewriting this relation in terms of
quantum minors we come to
\eqref{1rel}.
\epf

The next proposition is proved in the same way as its Yangian counterpart
\cite{nt:ry}; see also \cite{mno:yc}.

\bpr\label{prop:delqm}
The images of the quantum minors under the coproduct
are given by
\beql{delqm}
\Delta ({t\ts}^{a_1\cdots\ts a_r}_{b_1\cdots\ts b_r}(u))
=\sum_{c_1<\cdots < c_r}{t\ts}^{a_1\cdots\ts a_r}_{c_1\cdots\ts c_r}(u)
\ot {t\ts}^{c_1\cdots\ts c_r}_{b_1\cdots\ts b_r}(u),
\non
\end{equation}
where the summation is
over all subsets of indices $\{c_1,\dots,c_r\}$ from $\{1,\dots,n\}$.
\qed
\epr

\section{Gelfand--Tsetlin basis in $L(\la)$}
\label{sec:gt}

It has been observed in \cite{m:gt} that the raising and lowering operators
for the reduction $\gl_n\downarrow \gl_{n-1}$ (and more generally
for the corresponding Yangian reduction) can be given by
quantum minor formulas. This was used to construct
analogs of the Gelfand--Tsetlin basis for generic
Yangian modules. Here we modify
the arguments of \cite{m:gt} to
give $q$-analogs of the quantum minor formulas and construct
a basis of Gelfand--Tsetlin-type for the $\Uq$-module $L(\la)$.
Some other constructions of such bases
can be found in \cite{j:qr},
\cite{nt:yg}, \cite{t:ep}, \cite{uts:gz}.

A {\it pattern\/} $\La$ (associated with $\la$) is a sequence
of rows of integers $\La_n,\La_{n-1},\dots,\La_1$, where
$\La_r=(\lambda_{r1},\dots,\lambda_{rr})$ is the $r$-th row from the bottom,
the top row $\La_n$
coincides with $\la$, and the following
{\it betweenness conditions\/} are satisfied: for $r=1,\dots,n-1$
\beql{betw}
\la_{r+1,i+1}\leq\la_{ri}\leq \la_{r+1,i}
\qquad{\rm for}\quad i=1,\dots,r.
\end{equation}
We shall be using the notation $l_{ki}=\la_{ki}-i+1$.
Also, for any integer $m$ we set
\beql{int}
[m]=\frac{q^m-q^{-m}}{q-\qin}.
\end{equation}

\bpr\label{prop:gt}
There exists a basis $\{\xiL\}$ in $L(\la)$ parameterized
by the patterns $\La$ such that the action of the generators of
$\Uq$ is given by
\begin{align}\label{tgt}
t_k\ts\xiL&=q^{w_k}\ts\xiL,\qquad w_k=\sum_{i=1}^k\la_{ki}-\sum_{i=1}^{k-1}\la_{k-1,i},
\\
\label{egt}
e_k\ts\xiL&=-\sum_{j=1}^k\frac{[\ts l_{k+1,1}-l_{kj}]\cdots [\ts l_{k+1,k+1}-l_{kj}]}
{[\ts l_{k1}-l_{kj}]\cdots\wedge_j\cdots [\ts l_{kk}-l_{kj}]}
\ts \xi^{}_{\La+\delta_{kj}},\\
\label{fgt}
f_k\ts\xiL&=\sum_{j=1}^k\frac{[\ts l_{k-1,1}-l_{kj}]\cdots [\ts l_{k-1,k-1}-l_{kj}]}
{[\ts l_{k1}-l_{kj}]\cdots\wedge_j\cdots [\ts l_{kk}-l_{kj}]}
\ts \xi^{}_{\La-\delta_{kj}},
\end{align}
where $\La\pm\delta_{kj}$ is obtained from
$\La$ by replacing the entry $\lambda_{kj}$ with $\lambda_{kj}\pm1$,
and $\xi^{}_{\La}$ is supposed to be equal to zero
if $\La$ is not a pattern; the symbol $\wedge_j$ indicates
that the $j$-th factor is skipped.
\epr

\bpf Our proof of this result is based on the relations between
the quantum minors given in Section~\ref{sec:qmr}.
Set
\beql{ttil}
T_{ij}(u)=\frac{u\ts t_{ij}-u^{-1}\ts \bar t_{ij}}{q-\qin}.
\end{equation}
Clearly, $(q-\qin)\ts T_{ij}(u)=u\ts\ev(t_{ij}(u^2))$;
see \eqref{eval1}.
We also define the corresponding
quantum minors
${\ti\ts}^{a_1\cdots\ts a_r}_{b_1\cdots\ts b_r}(u)$ by the formulas
\eqref{qminor} or \eqref{qminor2} where all series $t_{ij}(u)$
are respectively replaced by $\ti_{ij}(u)$. Now for any
$1\leq a<r\leq n$
introduce the {\it lowering operators\/} by
\beql{lower}
\tau^{}_{ra}(u)=q^{r-a}\ts
{\ti\ts}^{a+1\ts\cdots\ts r}_{a\ts\cdots\ts r-1}(u).
\end{equation}
Note that by \eqref{center} we have
\beql{commu}
\tau^{}_{ra}(u)\ts\tau^{}_{sb}(v)=\tau^{}_{sb}(v)\ts\tau^{}_{ra}(u),
\end{equation}
if $b\leq a$ and $s\geq r$.

Let $\mu=(\mu_1,\dots,\mu_{n-1})$ be an $(n-1)$-tuple of integers
satisfying the inequalities
\beql{ineqmu}
\la_{i+1}\leq\mu_i\leq\la_{i},\qquad i=1,\dots,n-1.
\end{equation}
Introduce the vector
\beql{ximu}
\xi_{\mu}=\prod_{a=1}^{n-1}\tau_{na}(q^{-\mu_a-1})\cdots
\tau_{na}(q^{-\la_a+1})\ts \tau_{na}(q^{-\la_a})\ts\xi,
\end{equation}
where $\xi$ is the highest vector of $L(\la)$. An easy induction
with the use of \eqref{aq=bq} shows that $\xi_{\mu}$ satisfies
\beql{hvsub}
e_k\ts \xi_{\mu}=0,\quad k=1,\dots,n-2\qquad\text{and}\qquad
t_k\ts \xi_{\mu}=q^{\ts\mu_k} \ts \xi_{\mu},\quad k=1,\dots,n-1.
\end{equation}
Thus, if $\xi_{\mu}$ is nonzero then it generates a
$\U_q(\gl_{n-1})$-submodule isomorphic to $L(\mu)$.
Now given a pattern $\La$, we define vectors $\xiL\in L(\la)$ by
\beql{xildef}
\xiL=\prod_{r=2,\dots,n}^{\rightarrow}\prod_{a=1}^{r-1}
\tau^{}_{ra}(q^{-\la_{r-1,a}-1})\cdots \tau^{}_{ra}(q^{-\la_{ra}+1})
\ts \tau^{}_{ra}(q^{-\la_{ra}})\ts \xi.
\end{equation}
The relation \eqref{tgt} easily follows from \eqref{hvsub}
and the defining relations in $\Uq$.
We now derive the formulas \eqref{egt} and \eqref{fgt}
which together with \eqref{tgt} will imply that the vectors $\xiL$
are linearly independent and form a nontrivial submodule of $L(\la)$.
Since $L(\la)$ is irreducible this submodule must coincide with
$L(\la)$.
Below we shall only give a derivation of \eqref{egt}; the proof of
\eqref{fgt} is quite similar
and will be omitted. Note first, that, as follows from
\eqref{aq=bq}, if $k\geq r$ then $e_k$ commutes with
the lowering operator $\tau_{ra}(u)$. Therefore, we only need
to apply $e_{n-1}$ to the vector $\xi_{\mu}$ defined in
\eqref{ximu}. Let $a\in\{1,\dots,n-1\}$ be the least index
such that $\la_a-\mu_a>0$. We use a reverse induction on $a$
with the trivial base $a=n$ (i.e. $\xi_{\mu}=\xi$).
For a nonnegative integer $m$,
we introduce the products of the lowering operators by
\beql{Tau}
\Tau_{na}(u,m)=\prod_{i=1}^m \tau_{na}(q^{i-1}\ts u).
\end{equation}
We then have
\beql{xia}
\xi_{\mu}=\Tau_{na}(q^{-\la_a},\la_a-\mu_a)\ts \xi_{\mu'},
\end{equation}
where $\mu'$ is obtained from $\mu$ by replacing $\mu_a$ with $\la_a$.
We derive from \eqref{aq=bq} that
\beql{etau}
e_{n-1}\ts \tau_{na}(u)=\qin\ts \tau_{na}(u)\ts e_{n-1}
-q^{n-a-1}\ts {\ti\ts}^{a+1\ts\cdots\ts n}_{a\ts\cdots\ts n-2,n}(u).
\end{equation}
Then by induction we obtain
\beql{eTau}
\bal
e_{n-1}\ts \Tau_{na}(u,m)&=q^{-m}\ts \Tau_{na}(u,m)\ts e_{n-1}\\
{}&-\sum_{i=1}^m q^{n-a-i}\ts \tau_{na}(u)\cdots
{\ti\ts}^{a+1\ts\cdots\ts n}_{a\ts\cdots\ts n-2,n}(q^{i-1}\ts u)
\cdots \tau_{na}(q^{m-1}\ts u).
\eal
\end{equation}
Consider now the subalgebra $\Y_a$ of $\Yq$ generated by the coefficients
of $t_{ij}(u)$
with $a\leq i,j\leq n$. Using the relations \eqref{R21} for this subalgebra
we get
\beql{hata}
{\ti\ts}^{a+1\ts\cdots\ts n}_{a\ts\cdots\ts n-2,n}(u)\ts\tau_{na}(q\ts u)
=q\ts \tau_{na}(u)\ts {\ti\ts}^{a+1\ts\cdots\ts n}_{a\ts\cdots\ts n-2,n}(q\ts u).
\end{equation}
This brings \eqref{eTau} to the form
\beql{eTau2}
\bal
e_{n-1}\ts \Tau_{na}(u,m){}&
=q^{-m}\ts \Tau_{na}(u,m)\ts e_{n-1}\\
{}&-[m]\ts q^{n-a-1}\ts \Tau_{na}(u,m-1)\ts
{\ti\ts}^{a+1\ts\cdots\ts n}_{a\ts\cdots\ts n-2,n}(q^{m-1}\ts u).
\eal
\end{equation}
By the induction hypothesis, the action of $e_{n-1}$ on $\xi_{\mu'}$
is found from \eqref{egt}. So we only need to
calculate
${\ti\ts}^{a+1\ts\cdots\ts n}_{a\ts\cdots\ts n-2,n}(u)\ts    \xi_{\mu'}$
at $u=q^{-\mu_a-1}$. For this we use Proposition~\ref{prop:1rel}.
Clearly, the relation \eqref{1rel} remains valid if we replace each
quantum minor with the corresponding minor in the elements $\ti_{ij}(u)$.
By the induction hypothesis, we have
\beql{actmin}
{\ti\ts}^{a+1\ts\cdots\ts n-1}_{a+1\ts\cdots\ts n-1}(q^{-\mu_a-1})\ts   \xi_{\mu'}
=\prod_{i=a+1}^{n-1}[\ts m_i-m_a]\ts\xi_{\mu'},
\end{equation}
where $m_i=\mu_i-i+1$. Similarly, since
${\ti\ts}^{a+1\ts\cdots\ts n}_{a+1\ts\cdots\ts n}(u)$ commutes with the lowering
operators $\tau_{nb}(v)$, we obtain
\beql{actmin2}
{\ti\ts}^{a+1\ts\cdots\ts n}_{a+1\ts\cdots\ts n}(q^{-\mu_a-1})\ts   \xi_{\mu'}
=\prod_{i=a+1}^{n}[\ts l_i-m_a]\ts\xi_{\mu'}.
\end{equation}
Moreover, by \eqref{aq=bq} we have
\beql{qcom}
{\ti\ts}^{a+1\ts\cdots\ts n-1}_{a+1\ts\cdots\ts n-2,n}(u)
=[{\ti\ts}^{a+1\ts\cdots\ts n-1}_{a+1\ts\cdots\ts n-1}(u),e_{n-1}]_q,
\end{equation}
and so, the action of
${\ti\ts}^{a+1\ts\cdots\ts n-1}_{a+1\ts\cdots\ts n-2,n}(q^{-\mu_a-1})$
on $\xi_{\mu'}$ is also found by induction with the use of \eqref{actmin}.
It is now a matter of a straightforward calculation to check that the
resulting expression for the matrix elements of $e_{n-1}$ agrees with
\eqref{egt}. \epf

\section{Proof of Theorem~\ref{thm:main}}
\label{sec:proof}

Here we outline the main arguments in the proof of Theorem~\ref{thm:main}.
They closely follow the proof of its Yangian version of \cite{m:ic}
with the use of the quantum minor
relations given in Section~\ref{sec:qmr}.
For a pair of indices $i<j$ we shall denote
\beql{lanran}
\bal
\lan l_j,l_i\ran&=\{l_j,l_j+1,\dots,l_i\}\setminus\{l_j,l_{j-1},\dots,l_i\},
\\
\lan m_j,m_i\ran&=\{m_j,m_j+1,\dots,m_i\}\setminus\{m_j,m_{j-1},\dots,m_i\},
\eal
\end{equation}
where $l_i=\la_i-i+1$ and $m_i=\mu_i-i+1$.
In particular, if $\la_i=\la_{i+1}=\cdots=\la_j$ then
$\lan l_j,l_i\ran =\emptyset$.
It was shown in \cite[Proposition~2.8]{m:ic}
that the condition of Theorem~\ref{thm:main}
is equivalent to the following:
for all pairs of indices $1\leq i< j\leq n$
we have
\beql{ijcondeq}
m_j,m_i\not\in\lan l_j,l_i\ran\qquad\text{or}
\qquad l_j,l_i\not\in\lan m_j,m_i\ran.
\end{equation}
We start by proving that these conditions are sufficient
for the irreducibility of the $\Yq$-module $L(\la)\ot L(\mu)$.
Let $\xi$ and $\xi'$ denote the highest vectors of the $\gl_n$-modules
$L(\la)$ and $L(\mu)$, respectively.
The key part of the proof of sufficiency of the conditions is to show by induction
on $n$ that if $\zeta\in L(\la)\ot L(\mu)$ is a nonzero vector
satisfying \eqref{yanhv} for some series $\nu_i(u)$ then
\beql{sing}
\zeta={\rm const}\cdot\xi\ot\xi'.
\end{equation}
Then considering dual modules we also show that
the vector $\xi\ot\xi'$ is cyclic.

Note that the modules $L(\la)\ot L(\mu)$ and $L(\mu)\ot L(\la)$
are simultaneously reducible or irreducible. This can be easily deduced
from the formulas \eqref{hvtpr} for the highest weight
of the irreducible module \eqref{ten2}. So, we may assume
without loss of generality that
\beql{1ncond}
m_1,m_n\not\in\lan l_n,l_1\ran.
\end{equation}
Consider
the Gelfand--Tsetlin basis $\{\xiL\}$ of the
$\Uq$-module $L(\la)$; see Section~\ref{sec:gt}. The singular vector
$\zeta$ is uniquely written in the form
\beql{zeta}
\zeta=\sum_{\La}\xiL\ot \etaL,
\end{equation}
summed over all patterns $\La$ associated with $\la$, and
$\etaL\in L(\mu)$. We define the weight of a pattern $\La$
as the $n$-tuple $w(\La)=(w_1,\dots,w_n)$ where the $w_k$ are given in
\eqref{tgt}. We use a standard partial ordering on the weights
such that $w\preceq w'$ if and only if $w'-w$ is a $\ZZ_+$-linear
combination of the elements $\ve_i-\ve_{i+1}$, where $\ve_i$ is the $n$-tuple
with $1$ on the $i$-th place and zeros elsewhere.
Choose a minimal pattern $\La^0$ with respect to this ordering
among those occurring in the expansion \eqref{zeta}.
Then, exactly as in \cite[Lemmas~3.2 \& 3.3]{m:ic},
we show that $\etaLo$
is proportional to the highest vector $\xi'$ of $L(\mu)$
and that $\La^0$ is determined uniquely. Furthermore,
we apply Proposition~\ref{prop:gt} to demostrate that
due to the condition \eqref{1ncond}, the $(n-1)$-th row of $\La^0$
is ${\la_-{:=}}\tss(\la_1,\dots,\la_{n-1})$. This means that
the vector \eqref{zeta} belongs to
the $\Y_q(\gl_{n-1})$-span
of the vector $\xi\ot\xi'$, which is isomorphic to the tensor product
$L(\la_-)\ot L(\mu_-)$. Since the conditions \eqref{ijcondeq}
hold for $\la_-$ and $\mu_-$, we conclude by induction that \eqref{sing}
holds.

The next step is to show that under the conditions
\eqref{ijcondeq} (assuming \eqref{1ncond}   as well)
the vector $\xi\ot\xi'$ of the $\Yq$-module
$L=L(\la)\ot L(\mu)$ is cyclic.
The cyclicity of $L$ is equivalent
to the cocyclicity of the dual module $L^*= L(\la)^*\ot L(\mu)^*$
(that is, to the fact that any singular vector of $L^*$
is proportional to $\xi^*\ot\xi^{\prime\ts *}$).
To define the dual modules $L(\la)^*$ and $L(\mu)^*$
we use the anti-automorphism $\sigma$ of
$\Uq$ defined by
\beql{sigma}
\sigma :e_i\mapsto -e_i,\qquad \sigma :f_i\mapsto -f_i,\qquad
\sigma :t_i\mapsto t_i^{-1}.
\end{equation}
The dual space $L(\la)^*$
becomes a $\Uq$-module if we set
\beql{defduu}
(yf)(v)=f(\sigma(y)v),\qquad y\in\Uq,\quad f\in L(\la)^*,\quad v\in L(\la).
\end{equation}
It is easy to see that
the $\Uq$-module $L(\la)^*$ is isomorphic to
$L(-\la^{\omega})$, where $\la^{\omega}=(\la_n,\dots,\la_1)$, and so
\beql{isolst}
L^*\simeq L(-\la^{\omega})\ot L(-\mu^{\omega}).
\end{equation}
Next we verify that if $N$ is any submodule of $L$ then its
annihilator
\beq
\text{\rm Ann\ts} N=\{f\in L^*\ |\ f(v)=0\quad\text{\rm for all}\quad v\in N\}
\non
\end{equation}
is a nonzero submodule in $L^*$.
The claim now follows from the fact
that
if $\zeta'$ is a lowest singular vector of the module \eqref{isolst}
then $\zeta'$ is proportional to $\eta\ot\eta'$,
where $\eta$ and $\eta'$ are the lowest vectors of $L(\la)$ and $L(\mu)$,
respectively. This is proved by repeating the above argument
for the singular vector $\zeta$.

To prove the necessity of the conditions of the theorem
we use induction on $n$ again. It is not difficult to see that
if the $\Yq$-module $L(\la)\ot L(\mu)$ is irreducible then so are
the $\Y_q(\gl_{n-1})$-modules $L(\la_1,\dots,\la_{n-1})\ot L(\mu_1,\dots,\mu_{n-1})$
and $L(\la_2,\dots,\la_{n})\ot L(\mu_2,\dots,\mu_{n})$. Therefore, by the induction
hypothesis, the conditions \eqref{ijcondeq} can only be violated for
$i=1$ and $j=n$. Suppose this is the case.
Swapping $\la$ and $\mu$ if necessary, we may assume that
$m_n\in\lan l_n,l_1\ran$ and
$l_1\in\lan m_n,m_1 \ran$. There are two cases. First,
\beql{pcond}
m_n\in\lan l_{n},l_{n-1}\ran\qquad\text{and}\qquad l_1\in\lan m_{2},m_1\ran.
\end{equation}
Then there exist indices $r$ and $s$ such that
\beql{indrs}
m_2,\dots,m_r\in\{l_2,\dots,l_s\},\qquad
l_{s+1},\dots,l_{n-1}\in\{m_{r+1},\dots,m_{n-1}\}.
\non
\end{equation}
In particular, this implies that
\beql{ineqml}
l_i-m_i\in\ZZ_+ \qquad\text{for all}\quad i=2,\dots,n-1.
\end{equation}
The second case is
\beql{pcond2}
m_n\in\lan l_{p+1},l_p\ran\qquad\text{and}\qquad l_1\in\lan m_{n-p+1},m_{n-p}\ran
\end{equation}
for some $1\leq p\leq n-2$.
Then
$
l_{p-i+1}=m_{n-i}
$
for $i=1,\dots,p-1$.

We show that in both cases the module $L(\la)\ot L(\mu)$
contains a singular vector which is not proportional
to $\xi\ot\xi'$. Indeed, recall that $t_{ij}(u)$ acts
in the tensor product as the operator \eqref{tenoper} with $a=a'=1$.
Therefore, the operator $T_{ij}(u)=u^2\ts t_{ij}(u^2)$ in $L(\la)\ot L(\mu)$
is polynomial in $u$.
By analogy with \eqref{lower}, we introduce the lowering operators
\beql{lower2}
\tau^{}_{ra}(u)=
{T\ts}^{a+1\ts\cdots\ts r}_{a\ts\cdots\ts r-1}(u)
\end{equation}
and their products
\beql{Tau2}
\Tau_{na}(u,k)=\prod_{i=1}^k \tau_{na}(q^{i-1}\ts u),
\end{equation}
where $k$ is a nonnegative integer.
The numbers
\beql{ki}
k_i=l_i-m_{n-p+i},\qquad i=1,\dots,p
\non
\end{equation}
are positive integers and we define the vector $\theta\in L(\la)\ot L(\mu)$ by
\beql{thetot}
\theta=\Tau^{}_{n-p+1,1}(q^{-\la_1},k_1)\ts
\Tau^{\ts\prime}_{n-p+2,2}(q^{-\la_2},k_2)\cdots
\Tau^{\ts\prime}_{np}(q^{-\la_p},k_p)\ts(\xi\ot\xi'),
\non
\end{equation}
where $\Tau^{\ts\prime}_{na}(u,k)$
is the derivative of the polynomial
$\Tau_{na}(u,k)$. We need to prove that $\theta$ is annihilated by all operators
$T_{ij}(u)$ with $i<j$. It suffices to show that
\beql{mianth}
{T\ts}^{1\ts\cdots\ts k}_{1\ts\cdots\ts k-1,k+1}(u) \ts\theta =0,
\qquad k=1,\dots,n-1.
\end{equation}
Note that by \eqref{aq=bq} we have
\beql{becom}
{T\ts}^{1\ts\cdots\ts k}_{1\ts\cdots\ts k-1,k+1}(u)=
[{T\ts}^{1\ts\cdots\ts k}_{1\ts\cdots\ts k}(u),e_k]_q.
\end{equation}
Since the element ${t\ts}^{1\ts\cdots\ts k}_{1\ts\cdots\ts k}(u)$
is central in the subalgebra of $\Yq$ generated by $t_{ij}(u)$ with
$1\leq i,j\leq k$ the calculation is essentially reduced to showing
that $e_k\ts\theta=0$ for all $k$. The action of $e_k$ on
$\theta$ is found by a modified version of the argument which we used
in the derivation of \eqref{egt}; cf. \cite[Lemma~4.6]{m:ic}.

Finally, to prove that $\theta\ne 0$ we write it as a linear combination
\beql{gtLM}
\theta=\sum_{\La,M}c^{}_{\La,M}\ts\xiL\ot\xi'_M,
\end{equation}
where $\xiL$ and $\xi'_M$ are the Gelfand--Tsetlin basis vectors in $L(\la)$ and
$L(\mu)$. Applying Proposition~\ref{prop:delqm} we show that
\eqref{gtLM} has the form
\beql{gtLM2}
\theta=c\cdot \xiL\ot\xi'+\cdots,
\end{equation}
where $c$ is a nonzero constant and
\beql{xiLgt}
\xiL=\Tau^{}_{n-p+1,1}(q^{-\la_1},k_1)\ts
\Tau^{}_{n-p+2,2}(q^{-\la_2},k_2)\cdots
\Tau^{}_{np}(q^{-\la_p},k_p)\ts\xi
\end{equation}
is a vector of the Gelfand--Tsetlin basis of $L(\la)$; see \eqref{xildef}.
Thus, $\theta\ne 0$ which completes the proof of Theorem~\ref{thm:main}.

\end{document}